\documentclass{elsart}

\usepackage{graphics}
\usepackage{graphicx,float}
\usepackage{epsfig}
\usepackage{amsmath,amsfonts,amsgen}
\usepackage{dcolumn}
\usepackage[all]{xy}
\usepackage{makeidx,color}
\usepackage{color}

\newtheorem{definition}{Definition}
\newtheorem{property}{Property}

\newcolumntype{d}[1]{D{.}{\cdot}{#1}}
\newcolumntype{.}{D{.}{.}{-1}}
\DeclareMathOperator{\sinc}{sinc}

\newcommand\tabcaption{\def\@captype{table}\caption}
\newcommand{\be}{\begin{eqnarray*}}
\newcommand{\ee}{\end{eqnarray*}}
\newcommand{\ibe}{\begin{eqnarray}}
\newcommand{\iee}{\end{eqnarray}}

\def\3hf{\frac{3}{2}}

\def\np1{n+1}

\def\ip3hf{i+\frac{3}{2}}
\def\jp3hf{j+\frac{3}{2}}

\begin{document}
\begin{frontmatter}

\title{A new fractional derivative involving the normalized sinc function
without singular kernel}
\author{Xiao-Jun Yang$^{1,2}$, Feng Gao$^{1,2}$, J. A. Tenreiro Machado$^{3}$, Dumitru Baleanu$^{4,5}$}
\corauth[cor1]{Corresponding author:E-Mail: dyangxiaojun@163.com (X. J. Yang)}

\address{$^{1}$ School of Mechanics and Civil Engineering, China University of Mining and Technology, Xuzhou 221116, China}

\address{$^{2}$ State Key Laboratory for Geomechanics and Deep Underground Engineering,
China University of Mining and Technology, Xuzhou 221116, China}

\address{$^{3}$ Institute of Engineering, Polytechnic of Porto, Department of Electrical
Engineering, Rua Dr. Ant\'{o}nio Bernardino de Almeida, 4249-015 Porto, Portugal}

\address{$^{4}$ Department of Mathematics, Cankya University, Ogretmenler Cad. 14, Balgat-06530, Ankara Turkey}

\address{$^{5}$ Institute of Space Sciences, Magurele-Bucharest, Romania}

\begin{abstract}
In this paper, a new fractional derivative involving the normalized $\sinc$
function without singular kernel is proposed. The Laplace
transform is used to find the analytical solution of the anomalous
heat-diffusion problems. The comparative results between classical and
fractional-order operators are presented. The results are significant in the analysis
of one-dimensional anomalous heat-transfer problems.
\end{abstract}


\end{frontmatter}

\section{Introduction}

In recent years, fractional derivatives (FDs) in the sense of Caputo type
have used to describe anomalous behaviors of diffusive phenomena in
mathematical physics involving different kernels, such as the power-law
\cite{1}, exponential \cite{2}, Mittag-Leffler \cite{3}, stretched exponential \cite{4}, and
stretched Mittag-Leffler \cite{5} functions. For example, the fractional
diffusion-wave, in the power-law function kernel was considered in \cite{6}. The
numerical solution for the space-fractional diffusion equation was presented
in \cite{7}. The Cauchy problem for the time-fractional diffusion equations was
investigated in \cite{8}. With the help of the FD involving the
exponential-function kernel, the heat-diffusion problem with respect to a
non-singular fading memory was proposed in \cite{9}. The heat transfer problem
within the non-singular second grade fluid was discussed in \cite{10}. The
non-singular unsteady flow of the ordinary couple stress fluid was studied
in \cite{11}. With the use of the FD involving the stretched
Mittag-Leffler-function kernel, the Irving--Mullineux oscillator \cite{12} and the
Allen-Cahn equation \cite{13} were also analyzed. For more details see \cite{13a,13b}.

FDs in the sense of Riemann--Liouville type were developed in \cite{14,15,16}. We can
mention the studies about not only the Fokker-Planck \cite{17} and the diffusion \cite{18,19,20,21}
equations, but also the wave propagation \cite{22}. Furthermore, the Chen's system of the
Riemann--Liouville type \cite{23}, the~monotone iterative method
for neutral fractional differential equations \cite{24} and the
time-fractional-order Harry-Dym equation \cite{25} were also discussed.
Readers can find the more details about the distinct versions of FDs in \cite{26}.

The normalized $\sinc$ function, structured by Whittaker in \cite{27}, and its
properties were considered in \cite{28}. Furthermore, the Fourier \cite{29}, Laplace
\cite{29} and Sumudu \cite{30} transforms of the NSF were formulated. However, the FD involving
the normalized $\sinc$ function without singular kernel has not proposed.
Motivated by the idea, the present article derives a
new FD with respect to the normalized $\sinc$ function without singular kernel.
Furthermore, based on the new concept it is considered the applications in
one-dimensional anomalous heat-transfer problems.

The structure of the present paper is as follows. In Section 2, a new FD
with respect to the normalized $\sinc$ function without singular kernel is
presented. In Section 3, the anomalous heat-diffusion models and their
solutions are analyzed by means of the Laplace transform. Finally, the conclusion is
outlined in Section 4.

\section{Preliminaries, definitions and integral transforms}

In this section, we derive the FD involving the
normalized $\sinc$ function without singular kernel.

\subsection{A new FD involving the normalized $\sinc$
function without singular kernel}

\begin{definition}
The normalized $\sinc$ function is defined by \cite{27,28}:
\begin{equation}
\label{eq1}
\sinc\left( x \right)=\frac{\sin \left( {\pi x} \right)}{\pi x},
\end{equation}
where $x\in {\rm R}$.

If $\varphi \left( x \right)$ is any~smooth function with~compact support,
where $x\in {\rm R}$, then \cite{28}
\begin{equation}
\label{eq2}
\mathop {\lim }\limits_{\varpi \to 0} \frac{1}{\varpi }\sinc\left(
{\frac{x}{\varpi }} \right)=\mathop {\lim }\limits_{\varpi \to 0} \frac{\sin
\left( {\frac{\pi x}{\varpi }} \right)}{\pi x}=\delta \left( x \right),
\end{equation}
\begin{equation}
\label{eq3}
\sin c\left( 0 \right)=1,
\end{equation}
where
\begin{equation}
\label{eq4}
\mathop {\lim }\limits_{\varpi \to 0} \int\limits_{-\infty }^\infty
{\frac{\varphi \left( x \right)\sinc\left( {\frac{x}{\varpi }}
\right)}{\varpi }dx} =\varphi \left( 0 \right).
\end{equation}
\end{definition}

\begin{definition}
Let $\Pi \left( x \right)\in H^1\left( {a,b} \right)$ and $\mbox{ }b>a$. A
new FD involving the normalized $\sinc$ kernel of the function $\Pi \left( \mu \right)$ of order $\varpi
\mbox{ }\left( {\varpi \in \left( {0,1} \right)} \right)$ is defined as:
\begin{equation}
\label{eq5}
{ }_aD_\mu ^{\left( \varpi \right)} \Pi \left( \mu \right)=\frac{\varpi \wp
\left( \varpi \right)}{1-\varpi }\int\limits_a^\mu {\sinc\left(
{-\frac{\varpi \left( {\mu -x} \right)}{1-\varpi }} \right)\Pi ^{\left( 1
\right)}\left( x \right)dx} ,
\end{equation}
where $a\in \left( {-\infty ,\mu } \right)$, and $\wp \left( \varpi \right)$ is
a normalization constant depending on $\varpi $ such that $\wp \left( 0
\right)=\wp \left( 1 \right)=1$.
\end{definition}

Following Eq.(\ref{eq1}), we obtain
\begin{equation}
\label{eq6}
\mathop {\lim }\limits_{\varpi \to 0} \frac{1}{\varpi }\sinc\left(
{\frac{\mu -x}{\varpi }} \right)=\mathop {\lim }\limits_{\varpi \to 0}
\frac{\sin \left( {\frac{\pi \left( {\mu -x} \right)}{\varpi }} \right)}{\pi
\left( {\mu -x} \right)}=\delta \left( {\mu -x} \right),
\end{equation}
where $\varphi \left( x \right)$ is any~smooth function with~compact support
where $x\in {\rm R}_{ }$ such that
\begin{equation}
\label{eq7}
\mathop {\lim }\limits_{\varpi \to 0} \int\limits_{-\infty }^\infty {\varphi
\left( x \right)\frac{1}{\varpi }\sinc\left( {\frac{x-\mu }{\varpi }}
\right)dx} =\varphi \left( \mu \right).
\end{equation}
Thus, we have
\begin{equation}
\label{eq8}
\begin{array}{l}
 \mathop {\lim }\limits_{\varpi \to 0} { }_aD_\mu ^{\left( \varpi \right)}
\Pi \left( \mu \right)=\mathop {\lim }\limits_{\varpi \to 0} \frac{\varpi
\wp \left( \varpi \right)}{\left( {1-\varpi } \right)}\int\limits_a^\mu
{\sinc\left( {-\frac{\varpi \left( {\mu -x} \right)}{1-\varpi }} \right)\Pi
^{\left( 1 \right)}\left( x \right)dx} \\
 =\left( {\mathop {\lim }\limits_{\varpi \to 0} \wp \left( \varpi \right)}
\right)\int\limits_a^\mu {\delta \left( {\mu -x} \right)\Pi ^{\left( 1
\right)}\left( x \right)dx} \\
 =\Pi ^{\left( 1 \right)}\left( x \right). \\
 \end{array}
\end{equation}
When
\begin{equation}
\label{eq9}
\mathop {\lim }\limits_{\varpi \to 1} \frac{\varpi \wp \left( \varpi
\right)}{1-\varpi }\left[ {\sinc\left( {-\frac{\varpi \left( {\mu -x}
\right)}{1-\varpi }} \right)} \right]=\mathop {\lim }\limits_{\varpi \to
\mbox{1}} \wp \left( \varpi \right)\frac{\varpi }{1-\varpi }\frac{\sin
\left( {\frac{\pi \left( {\mu -x} \right)}{\frac{1-\varpi }{\varpi }}}
\right)}{\pi \left( {\mu -x} \right)}=\mathop {\lim }\limits_{\varpi \to
\mbox{1}} \mbox{1,}
\end{equation}
we have
\begin{equation}
\label{eq10}
\begin{array}{l}
 \mathop {\lim }\limits_{\varpi \to 1} { }_aD_\mu ^{\left( \varpi \right)}
\Pi \left( \mu \right)=\mathop {\lim }\limits_{\varpi \to 1} \frac{\varpi
\wp \left( \varpi \right)}{\left( {1-\varpi } \right)}\int\limits_a^\mu
{\sinc\left( {-\frac{\varpi \left( {\mu -x} \right)}{1-\varpi }} \right)\Pi
^{\left( 1 \right)}\left( x \right)dx} \\
 =\mathop {\lim }\limits_{\varpi \to 1} \int\limits_a^\mu {\Pi ^{\left( 1
\right)}\left( x \right)dx} \\
 =\Pi \left( \mu \right)-\Pi \left( a \right). \\
 \end{array}
\end{equation}
For $n\ge 1$ and $\varpi \in \left( {0,1} \right)$, the FD $D_\mu ^{\left( {n+\varpi } \right)} \Pi \left( \mu \right)$ of
order $\left( {n+\omega } \right)$ is defined as:
\begin{equation}
\label{eq11}
{ }_aD_\mu ^{\left( {n+\varpi } \right)} \Pi \left( \mu \right):={ }_aD_\mu
^{\left( n \right)} \left( {{ }_aD_\mu ^{\left( \varpi \right)} \Pi \left(
\mu \right)} \right).
\end{equation}

\begin{property}

(T1) ${ }_0D_\mu ^{\left( \varpi \right)} \theta =0$, where $\theta $ is a
constant;

(T2) ${ }_0D_\mu ^{\left( \varpi \right)} \mu =\frac{\varpi \wp \left(
\varpi \right)}{1-\varpi }\int\limits_0^\mu {\sinc\left( {\frac{\varpi
x}{1-\varpi }} \right)dx} $.
\end{property}

\begin{proof}
We have from Eq.(\ref{eq5}) that
\begin{equation}
\label{eq12}
{ }_0D_\mu ^{\left( \varpi \right)} \theta =\frac{\varpi \wp \left( \varpi
\right)}{1-\varpi }\int\limits_0^\mu {\sinc\left( {-\frac{\varpi \left(
{\mu -x} \right)}{1-\varpi }} \right)\theta ^{\left( 1 \right)}dx} =0.
\end{equation}
We have, by using the definition Eq.(\ref{eq5}),
\begin{equation}
\label{eq13}
{ }_0D_\mu ^{\left( \varpi \right)} \mu =\frac{\varpi \wp \left( \varpi
\right)}{1-\varpi }\int\limits_0^\mu {\sinc\left( {-\frac{\varpi \left(
{\mu -x} \right)}{1-\varpi }} \right)dx} =\frac{\varpi \wp \left( \varpi
\right)}{1-\varpi }\int\limits_0^\mu {\sinc\left( {\frac{\varpi x}{1-\varpi
}} \right)dx} .
\end{equation}
\end{proof}

\subsection{Integral transforms of the new FD
involving the normalized $\sinc$ function without singular kernel}

Here, we have \cite{29}
\begin{equation}
\label{eq14}
\aleph \left\{ {\sinc\left( x \right)} \right\}=\aleph \left\{ {\frac{\sin
\left( {\pi x} \right)}{\pi x}} \right\}=\sqrt {\frac{1}{2\pi }} H\left(
{\pi -\left| \xi \right|} \right)
\end{equation}
such that
\begin{equation}
\label{eq15}
\begin{array}{l}
 \aleph \left\{ {\sinc\left( {-\frac{\varpi x}{1-\varpi }} \right)}
\right\}=\aleph \left\{ {\frac{\sin \left( {-\frac{\varpi \pi x}{1-\varpi }}
\right)}{-\frac{\varpi \pi x}{1-\varpi }}} \right\} \\
 =-\sqrt {\frac{1}{2\pi }} \frac{1-\varpi }{\varpi }H\left( {-\frac{\varpi
\pi }{1-\varpi }-\left| \xi \right|} \right) \\
 =\sqrt {\frac{1}{2\pi }} \frac{1-\varpi }{\varpi }H\left( {\frac{\varpi \pi
}{1-\varpi }+\left| \xi \right|} \right), \\
 \end{array}
\end{equation}
where $\aleph $ is the Fourier transform operator \cite{29}, and $H\left( x
\right)$ is the Heaviside function \cite{29}.

The Fourier transform of Eq.(\ref{eq5}) can be written as
\begin{equation}
\label{eq16}
\begin{array}{l}
 \aleph \left\{ {{ }_0D_\mu ^{\left( \varpi \right)} \Pi \left( \mu \right)}
\right\} \\
 =\aleph \left\{ {\frac{\varpi \wp \left( \varpi \right)}{1-\varpi
}\int\limits_0^\mu {\sinc\left( {-\frac{\varpi \left( {\mu -x}
\right)}{1-\varpi }} \right)\Pi ^{\left( 1 \right)}\left( x \right)dx} }
\right\} \\
 =\frac{\varpi \wp \left( \varpi \right)}{1-\varpi }\aleph \left\{ {\sinc\left( {-\frac{\varpi x}{1-\varpi }} \right)} \right\}\aleph \left\{ {\Pi
^{\left( 1 \right)}\left( x \right)} \right\} \\
 =\frac{\varpi \wp \left( \varpi \right)}{1-\varpi }\left[ {\sqrt
{\frac{1}{2\pi }} \frac{1-\varpi }{\varpi }H\left( {\frac{\varpi \pi
}{1-\varpi }+\left| \xi \right|} \right)} \right]\left[ {i\xi \Pi \left( \xi
\right)} \right] \\
 =i\xi \sqrt {\frac{1}{2\pi }} \wp \left( \varpi \right)H\left(
{\frac{\varpi \pi }{1-\varpi }+\left| \xi \right|} \right)\Pi \left( \xi
\right), \\
 \end{array}
\end{equation}
where $\aleph \left\{ {\Pi \left( \mu \right)} \right\}=\Pi \left( \xi
\right)$.

Similarly, we have \cite{29}
\begin{equation}
\label{eq17}
\Im \left\{ {\sinc\left( x \right)} \right\}=\Im \left\{ {\frac{\sin \left(
{\pi x} \right)}{\pi x}} \right\}=\frac{1}{\pi }\mbox{tan}^{-1}\left(
{\frac{\pi }{s}} \right)
\end{equation}
such that
\begin{equation}
\label{eq18}
\Im \left\{ {\sinc\left( {-\frac{\varpi }{1-\varpi }x} \right)}
\right\}=\Im \left\{ {\frac{\sin \left( {-\frac{\varpi }{1-\varpi }\pi x}
\right)}{-\frac{\varpi \pi x}{1-\varpi }}} \right\}=\frac{1}{\frac{\varpi
\pi }{1-\varpi }}\mbox{tan}^{-1}\left( {\frac{\frac{\varpi \pi }{1-\varpi
}}{s}} \right),
\end{equation}
where $\Im $ is the Laplace transform operator \cite{29}.

From Eq.(\ref{eq18}) the Laplace transform of Eq.(\ref{eq5}) can be given by:
\begin{equation}
\label{eq19}
\begin{array}{l}
 \Im \left\{ {{ }_0D_\mu ^{\left( \varpi \right)} \Pi \left( \mu \right)}
\right\} \\
 =\Im \left\{ {\frac{\varpi \wp \left( \varpi \right)}{1-\varpi
}\int\limits_0^\mu {\sinc\left( {-\frac{\varpi \left( {\mu -x}
\right)}{1-\varpi }} \right)\Pi ^{\left( 1 \right)}\left( x \right)dx} }
\right\} \\
 =\frac{\varpi \wp \left( \varpi \right)}{1-\varpi }\Im \left\{ {\sinc\left( {-\frac{\varpi x}{1-\varpi }} \right)} \right\}\Im \left\{ {\Pi
^{\left( 1 \right)}\left( x \right)} \right\} \\
 =\frac{\wp \left( \varpi \right)}{\pi }\mbox{tan}^{-1}\left(
{\frac{\frac{\varpi \pi }{1-\varpi }}{s}} \right)\left( {s\Pi \left( s
\right)-\Pi \left( 0 \right)} \right), \\
 \end{array}
\end{equation}
where $\Im \left\{ {\Pi \left( \mu \right)} \right\}=\Pi \left( s \right)$.

As a direct result, we have \cite{30}
\begin{equation}
\label{eq20}
\Re \left\{ {\sinc\left( x \right)} \right\}=\Re \left\{ {\frac{\sin \left(
{\pi x} \right)}{\pi x}} \right\}=\frac{\mbox{tan}^{-1}\left( {\pi \zeta }
\right)}{\pi \zeta }
\end{equation}
such that
\begin{equation}
\label{eq21}
\Re \left\{ {\sinc\left( {-\frac{\varpi }{1-\varpi }x} \right)}
\right\}=\Im \left\{ {\frac{\sin \left( {-\frac{\varpi }{1-\varpi }\pi x}
\right)}{-\frac{\varpi \pi x}{1-\varpi }}}
\right\}=\frac{\mbox{tan}^{-1}\left( {\frac{\varpi \pi \zeta }{1-\varpi }}
\right)}{\frac{\varpi \pi \zeta }{1-\varpi }}_{,}
\end{equation}
where $\Re $ is the Sumudu transform operator \cite{30}.

Thus, we have from Eq.(\ref{eq13}) that
\begin{equation}
\label{eq22}
\begin{array}{l}
 \Re \left\{ {{ }_0D_\mu ^{\left( \varpi \right)} \Pi \left( \mu \right)}
\right\} \\
 =\Re \left\{ {\frac{\varpi \wp \left( \varpi \right)}{1-\varpi
}\int\limits_0^\mu {\sinc\left( {-\frac{\varpi \left( {\mu -x}
\right)}{1-\varpi }} \right)\Pi ^{\left( 1 \right)}\left( x \right)dx} }
\right\} \\
 =\frac{\varpi \wp \left( \varpi \right)}{1-\varpi }\Re \left\{ {\sinc\left( {-\frac{\varpi x}{1-\varpi }} \right)} \right\}\Re \left\{ {\Pi
^{\left( 1 \right)}\left( x \right)} \right\} \\
 =\frac{\wp \left( \varpi \right)}{\pi \zeta }\mbox{tan}^{-1}\left(
{\frac{\varpi \pi \zeta }{1-\varpi }} \right)\left( {\frac{\Pi \left( \zeta
\right)-\Pi \left( 0 \right)}{\zeta }} \right), \\
 \end{array}
\end{equation}
where $\Re \left\{ {\Pi \left( \mu \right)} \right\}=\Pi \left( \zeta
\right)$.

\section{Modelling the anomalous heat-diffusion problems }

In this section, we model the anomalous heat-diffusion problems involving
fractional-time and -space derivatives of the normalized $\sinc$ function
without singular kernel.

\textbf{Example 1}

The anomalous heat-diffusion within the fractional-time derivative of the
normalized $\sinc$ function without singular kernel is written as:
\begin{equation}
\label{eq23}
{ }_0D_\tau ^{\left( \varpi \right)} \Pi \left( {\mu ,\tau } \right)=\kappa
\frac{\partial ^2\Pi \left( {\mu ,\tau } \right)}{\partial \mu ^2},\mbox{
}\mu >0,\mbox{ }\tau >0,
\end{equation}
subjected to the initial and boundary conditions:
\begin{equation}
\label{eq24}
\Pi \left( {\mu ,0} \right)=0,\mbox{ }\mu >0,
\end{equation}
\begin{equation}
\label{eq25}
\Pi \left( {0,\tau } \right)=\lambda \left( \tau \right)\mbox{, }\tau >0,
\end{equation}
\begin{equation}
\label{eq26}
\Pi \left( {\mu ,\tau } \right)\to 0,\mbox{ }as\mbox{ }\mu \to \infty
,\mbox{ }\tau >0,
\end{equation}
where $\kappa $ is the thermal diffusivity.

With the aid of Eq.(\ref{eq19}), Eq.(\ref{eq23}) can be transferred into
\begin{equation}
\label{eq27}
\frac{\wp \left( \varpi \right)}{\pi }\mbox{tan}^{-1}\left(
{\frac{\frac{\varpi \pi }{1-\varpi }}{s}} \right)\left( {s\Pi \left( {\mu
,s} \right)-\Pi \left( {\mu ,0} \right)} \right)=\kappa \frac{d^2\Pi \left(
{\mu ,s} \right)}{d\mu ^2}.
\end{equation}
From Eq.(\ref{eq24}) we have the following:
\begin{equation}
\label{eq28}
\frac{d^2\Pi \left( {\mu ,s} \right)}{d\mu ^2}=\frac{\wp \left( \varpi
\right)s}{\pi \kappa }\mbox{tan}^{-1}\left( {\frac{\frac{\varpi \pi
}{1-\varpi }}{s}} \right)\Pi \left( {\mu ,s} \right),
\end{equation}
which leads to
\begin{equation}
\label{eq29}
\Pi \left( {\mu ,s} \right)=\Omega _1 \exp \left( {-\mu \sqrt {\rm H} }
\right)+\Omega _2 \exp \left( {\mu \sqrt {\rm H} } \right),
\end{equation}
where $\Omega _1 $ and $\Omega _2 $ are two unknown constants and
\begin{equation}
\label{eq30}
{\rm H}=\frac{\wp \left( \varpi \right)s}{\pi \kappa }\mbox{tan}^{-1}\left(
{\frac{\frac{\varpi \pi }{1-\varpi }}{s}} \right).
\end{equation}
In view of Eq.(\ref{eq25}) and Eq.(\ref{eq26}), we have
\begin{equation}
\label{eq31}
\Omega _2 =0
\end{equation}
such that
\begin{equation}
\label{eq32}
\Pi \left( {\mu ,s} \right)=\lambda \left( s \right)\exp \left( {-\mu \sqrt
{\rm H} } \right),
\end{equation}
where $\Im \left\{ {\lambda \left( \mu \right)} \right\}=\lambda \left( s
\right)$.

Thus, the Laplace transform solution of Eq.(\ref{eq23}) is:
\begin{equation}
\label{eq33}
\Pi \left( {\mu ,s} \right)=\lambda \left( s \right)\exp \left( {-\sqrt
{\frac{\wp \left( \varpi \right)s}{\pi \kappa }\mbox{tan}^{-1}\left(
{\frac{\frac{\varpi \pi }{1-\varpi }}{s}} \right)} \mu } \right).
\end{equation}

\textbf{Example 2}

The anomalous heat-diffusion within the fractional-space derivative of the
normalized $\sinc$ function without singular kernel is
\begin{equation}
\label{eq34}
\frac{\partial \Pi \left( {\mu ,\tau } \right)}{\partial \tau }=\kappa {
}_0D_\mu ^{\left( 1 \right)} \left( {{ }_0D_\mu ^{\left( \varpi \right)} \Pi
\left( {\mu ,\tau } \right)} \right),\mbox{ }\mu >0,\mbox{ }\tau >0,
\end{equation}
with the initial and boundary conditions:
\begin{equation}
\label{eq35}
\Pi \left( {\mu ,0} \right)=0,\mbox{ }\mu >0,
\end{equation}
\begin{equation}
\label{eq36}
\Pi \left( {0,\tau } \right)=\lambda \left( \tau \right)\mbox{, }\tau >0,
\end{equation}
\begin{equation}
\label{eq37}
\Pi \left( {\mu ,0} \right)\to 0,\mbox{ }as\mbox{ }\mu \to \infty ,\mbox{
}\tau >0,
\end{equation}
where $\kappa $ is the thermal diffusivity, and
\begin{equation}
\label{eq38}
{ }_0D_\mu ^{\left( 1 \right)} \left( {{ }_0D_\mu ^{\left( \varpi \right)}
\Pi \left( {\mu ,\tau } \right)} \right)=\frac{\varpi \wp \left( \varpi
\right)}{1-\varpi }\frac{\partial }{\partial \mu }\int\limits_0^\mu {\sinc
\left( {-\frac{\varpi \left( {\mu -x} \right)}{1-\varpi }} \right)\Pi
^{\left( 1 \right)}\left( {x,\tau } \right)dx} .
\end{equation}
With the help of Eq.(\ref{eq19}) and Eq.(\ref{eq35}), Eq.(\ref{eq34}) can be written as:
\begin{equation}
\label{eq39}
{ }_0D_\mu ^{\left( 1 \right)} \left( {{ }_0D_\mu ^{\left( \varpi \right)}
\Pi \left( {\mu ,s} \right)} \right)=\frac{s}{\kappa }\Pi \left( {\mu ,s}
\right),
\end{equation}
where
\begin{equation}
\label{eq40}
{ }_0D_\mu ^{\left( 1 \right)} \left( {{ }_0D_\mu ^{\left( \varpi \right)}
\Pi \left( {\mu ,s} \right)} \right)=\frac{\varpi \wp \left( \varpi
\right)}{1-\varpi }\frac{\partial }{\partial \mu }\int\limits_0^\mu {\sinc
\left( {-\frac{\varpi \left( {\mu -x} \right)}{1-\varpi }} \right)\Pi
^{\left( 1 \right)}\left( {x,s} \right)dx} .
\end{equation}
By the integration of Eq.(\ref{eq39}) we have
\begin{equation}
\label{eq41}
{ }_0D_\mu ^{\left( \varpi \right)} \Pi \left( {\mu ,s} \right)=\frac{\varpi
\wp \left( \varpi \right)}{1-\varpi }\int\limits_0^\mu {\sinc\left(
{-\frac{\varpi \left( {\mu -x} \right)}{1-\varpi }} \right)\Pi ^{\left( 1
\right)}\left( {x,s} \right)dx} =\int\limits_0^\mu {\frac{s}{\kappa }\Pi
\left( {x,s} \right)dx} +\Theta ,
\end{equation}
where $\Theta $ is a constant.

By taking the Sumudu transform operator with $\mu $ and $\Theta =0$, we have
\begin{equation}
\label{eq42}
\frac{\varpi \wp \left( \varpi \right)}{1-\varpi }\int\limits_0^\mu {\sinc\left( {-\frac{\varpi \left( {\mu -x} \right)}{1-\varpi }} \right)\Pi
^{\left( 1 \right)}\left( {x,s} \right)dx} =\int\limits_0^\mu
{\frac{s}{\kappa }\Pi \left( {x,s} \right)dx} ,
\end{equation}
which implies that
\begin{equation}
\label{eq43}
\frac{\wp \left( \varpi \right)}{\pi \zeta ^2}\mbox{tan}^{-1}\left(
{\frac{\varpi \pi \zeta }{1-\varpi }} \right)\left( {\Pi \left( {\zeta ,s}
\right)-\Pi \left( {0,s} \right)} \right)=\frac{s\zeta }{\kappa }\Pi \left(
{\zeta ,s} \right).
\end{equation}
From Eq.(\ref{eq35}) and Eq.(\ref{eq43}), we have the following:
\begin{equation}
\label{eq44}
\frac{\wp \left( \varpi \right)}{\pi \zeta ^2}\mbox{tan}^{-1}\left(
{\frac{\varpi \pi \zeta }{1-\varpi }} \right)\left( {\Pi \left( {\zeta ,s}
\right)-\lambda \left( s \right)} \right)=\frac{s\zeta }{\kappa }\Pi \left(
{\zeta ,s} \right).
\end{equation}
Thus, we have
\begin{equation}
\label{eq45}
\Pi \left( {\zeta ,s} \right)=\frac{\frac{\wp \left( \varpi \right)\lambda
\left( s \right)}{\pi \zeta ^2}\mbox{tan}^{-1}\left( {\frac{\varpi \pi \zeta
}{1-\varpi }} \right)}{\frac{\wp \left( \varpi \right)}{\pi \zeta
^2}\mbox{tan}^{-1}\left( {\frac{\varpi \pi \zeta }{1-\varpi }}
\right)-\frac{s\zeta }{\kappa }}
\end{equation}
where $\Re \left\{ {\Pi \left( {\mu ,s} \right)} \right\}=\Pi \left( {\zeta
,s} \right)$ represents the Sumudu transform operator \cite{30}.

From Eq.(\ref{eq45}), the Laplace transform solution of Eq.(\ref{eq23}) is:
\begin{equation}
\label{eq46}
\Pi \left( {\mu ,s} \right)=\Re ^{-1}\left\{ {\frac{\frac{\wp \left( \varpi
\right)\lambda \left( s \right)}{\pi \zeta ^2}\mbox{tan}^{-1}\left(
{\frac{\varpi \pi \zeta }{1-\varpi }} \right)}{\frac{\wp \left( \varpi
\right)}{\pi \zeta ^2}\mbox{tan}^{-1}\left( {\frac{\varpi \pi \zeta
}{1-\varpi }} \right)-\frac{s\zeta }{\kappa }}} \right\},
\end{equation}
where $\Re ^{-1}\left\{ {\Pi \left( {\zeta ,s} \right)} \right\}=\Pi \left(
{\mu ,s} \right)$ represents the inverse Sumudu transform operator [30].

When $\varpi =0$, Eq.(\ref{eq33}) and Eq.(\ref{eq46}) become the Laplace transform solution
of the classical heat-diffusion equation \cite{31}:
\begin{equation}
\label{eq47}
\Pi \left( {\mu ,s} \right)=\lambda \left( s \right)\exp \left( {-\sqrt
{\frac{s}{\kappa }} \mu } \right).
\end{equation}
which is in agreement with the result in \cite{29}.

\section{Conclusions}

In the present study, we addressed a new FD in respect to
the normalized $\sinc$ function without singular kernel.
Moreover, the Fourier, Laplace and Sumudu transforms of the FD operator
and the Laplace--transform solutions of the anomalous heat-diffusion equations were
considered. The analytical solutions of the
classical and anomalous heat-diffusion equations in the form of the Laplace
transform were also compared. The new formulation may be used to support a new
perspective for describing the anomalous behaviors in mathematical
physics.

%%%%%%%%%%%%%%%%%%%%%%%%%%%%%%%%%%%%%%%%%%


\begin{thebibliography}{31}
\bibitem{27}
Whittaker, E. T. (1915). On the functions which are represented by the expansions of the interpolation-theory.~Proceedings of the Royal Society of Edinburgh,~35, 181-194.
\end{thebibliography}
\end{document}